\documentclass[a4paper]{amsart}
\usepackage{commandoos}

\begin{document}
	\author{A. Eggink}
	\title{Calculating the Local Ideal Class Monoid and Gekeler Ratios} 
	\begin{abstract}
	Let $A = \F_q[T]$, $\p \subset A$ prime, $f(x) \in A[x]$ irreducible and set $R = A[x]/f(x)$. Denote its completion by $R_\p$. The ideal class monoid $\ICM(R_\p)$ is the set of fractional $R_\p$ ideals modulo the principal $R_\p$ ideals. We provide an algorithm to compute $\ICM(R_\p)$. In the process we also get algorithms to compute the overorders and weak equivalence classes of $R_\p$. We then use the algorithms to compute the product of local Gekeler ratios $\prod_{\p \subset A} v_\p(f) = \prod_{\p \subset A} \lim_{n \ra \infty} \frac{|\{M \in \Mat_r(A/\p^n)\mid \charpoly(M)=f\}}{|\SL_r(A/\p^n)|/|\p|^{n(r-1)}}$. This provides part of an algorithm to compute the weighted size of an isogeny class of Drinfeld modules.  
	\end{abstract}
	
	\maketitle

\section{Introduction}

This paper consists of two parts. In the first part we describe an algorithm to compute the ideal class monoid of the completion at a prime $\p \subset \F_q[T]$ of the order $\F_q[T][x]/f(x)$ for some irreducible $f$. In the second part we provide an application of this to compute the product of local Gekeler ratios. 

\subsection{Ideal Class Monoid}
Write $A = \F_q[T]$ and take $\p \subset A$ prime and $f \in A[x]$ irreducible. Write $K$ for the fraction field of $A[x]/f(x)$ and let $R$ be an order in $K$. Denote with $R_\p$ its completion and with $R_{(\p)}$ its localisation at $\p$. The ideal class monoid of $S$, where $S$ is one of $R,R_\p$ or $R_\p$, is 
$$\ICM(S) = \{I \text{ fractional ideal of }S\}/\{J \text{ principal ideal of }S\}.$$ 
In \cite{Stefano} Marseglia describes, using \cite{PicGrpRef} and \cite{WeakEquivRef}, algorithms to compute $\ICM(R)$ in the number field setting. These algorithms also work in the function field setting, with one change in the inseparable case, see \th\ref{ICMWeakEquivInsep}. We will follow the same strategy for the local case. Let $I$ and $J$ be fractional ideals of $S$ where $S$ is one of $R,R_\p$ or $R_{(\p)}$. We say that $I$ and $J$ are weakly equivalent if for all primes $\q$ of $S$, there exists some $x \in K_\q$ such that $I_\q = xJ_\q$. We write $[I]$ for the weak equivalence class of $I$ and $(I:I)$ for the multiplicator ring of $I$. Let $S'$ be an overorder of $S$. We denote the set of weak equivalence classes in $S$ by $W(S)$ and set $W_{S'}(S) = \{[I] \in W(S)\mid (I:I) = S'\}$. In the global case we have by \cite[Theorem 4.6]{Stefano} that $\ICM(R) = \sqcup_{S} W_S(R)\Pic(S) = \sqcup_S \{ \{IJ\} \mid I \in W_S(R),\ J \in \Pic(S)\}$, where the union runs over the overorders of $S$. Since in the local setting the Picard groups are trivial, see \th\ref{spPicTrivRp}, we have the simplified expression $\ICM(R_\p) = \sqcup_S W_S(R_\p)$, where the union runs over all $A_\p$-overorders of $R_\p$. This strategy leads to the following algorithms:
\begin{enumerate}
\item Given $R = A[x]/f(x)$ with $f$ irreducible and a $\p \subset A$ prime, an algorithm to compute a set $\mathcal{S}$ of $A$-overorders of $R$ such that we have a bijection:
$$\mathcal{S} \ra \{A_\p\text{-overorders of }R_\p\},\ S \ra S_\p.$$ \label{alg1}
\item Given any $A$-order $R$, a prime $\p \subset A$ and an $A$-overorder $S$ of $R$, an algorithm to compute a set $\mathcal{I}$ of fractional $R$-ideals such that we have a bijection
$$\mathcal{I} \ra W_{S_\p}(R_\p),\ I \ra [I_\p].$$ \label{alg2} 
\item Given $R = A[x]/f(x)$ with $f$ irreducible and a prime $\p \subset A$, an algorithm to compute a set $\mathcal{I}$ of fractional $R$-ideals such that we have a bijection 
$$\mathcal{I} \ra \ICM(R_\p),\ I \ra \{I_\p\}.$$ \label{alg3} 
\end{enumerate}
Algorithms \eqref{alg1},\eqref{alg2} and \eqref{alg3} are explained in paragraphs \ref{parICOO}, \ref{parICWE} and \ref{parICICM} respectively. Paragraph \ref{parIClemmas} contains some preliminary results. 

\subsection{Gekeler Ratios}
We will use the algorithm to compute local ideal class monoids to compute the product of Gekeler ratios. The Gekeler ratios are related to counting Drinfeld modules. In order to explain this, we will introduce Drinfeld modules and some related concepts. Let $k$ be a field of characteristic $p$. The twisted polynomial ring is $k\{\t\} = \{ \sum_{i=0}^n a_i\t^i \mid n \in \N,\ a_i \in k\}$, with the usual addition, but multiplication defined by the rule $\t a = a^q \t$ for all $a\in k$. A Drinfeld module is an $\F_q$-algebra homomorphism $\phi: A \ra k\{\t\}$ such that $\phi(A) \not\subset k$. Note that $\phi$ is determined by $\phi(T) = \phi_T$. The $\tau$-degree of $\phi_T$ is called the rank of $\phi$. Rank two Drinfeld modules are the function field analogues of elliptic curves, while higher rank Drinfeld modules have similarities to higer dimensional abelian varieties. Just as in that case, we have isogeny and isomorphism classes of Drinfeld modules. We say that the Drinfeld modules $\phi$ and $\psi$ are isogenous if there exists a $u \in k\{\t\}\bs\{0\}$ such that $u\phi_T = \psi_T u$, and isomorphic if there exists a $c \in k\{\t\}^\ti = k^\ti$ such that $c\phi_T = \psi_Tc$. (Being isogenous is an equivalence relation, since there always exists a dual isogeny $\hat{u}$ such that $\phi_T\hat{u} = \hat{u}\psi_T$, see \cite[Proposition 3.3.12]{DM}.) We are interested in Drinfeld modules over some finite field, so $k = \F_{q^n}$. Since isogenies preserve rank, this implies that every isogeny class is finite. This leads to the question: How big is an isogeny class? We will look at the weighted seize of an isogeny class, which is $h^{*}(\phi) = \sum_{\psi} \frac{1}{|\Aut_{\F_{q^n}}(\psi)|}$, where the sum is taken over all isomorphism classes $\psi$ such that $\psi$ is isogenous to $\phi$. The theory of Drinfeld modules over finite fields tells us that the isogeny class of $\phi$ is determined by the characteristic polynomial of its Frobenius, see \cite[Theorem 4.3.2]{DM}. The characteristic polynomial is an element of $A[x]$ that can be computed in SageMath or Magma. For ordinary Drinfeld modules with a commutative endomorphism ring, Bray proved in \cite{ThesisAmie} that 
$$h^*(\phi) = v_{\infty}(f) \prod_{\p \subset A} v_\p(f).$$
Here $f$ is the characteristic polynomial of the Frobenius of $\phi$, $v_\infty(f)$ is a global term and 
$$v_\p(f) = \lim_{n \ra \infty} \frac{|\{M \in \Mat_r(A/\p^n)\mid \charpoly(M)=f\}}{|\SL_r(A/\p^n)|/|\p|^{n(r-1)}}$$ 
is the local Gekeler ratio at $\p$. Since $\phi$ has a commutative endomorphism ring, $f$ is irreducible. There did not exist algorithms to compute $v_\infty(f)$ or $v_\p(f)$ before, but now we have:
\begin{thm}
There exists an algorithm that given an irreducible $f \in A[x]$ computes $\prod_{\p \subset A} v_\p(f)$. 
\end{thm}
For the proof we follow ideas from \cite[Lemma 6.4 and Corrolary 6.5]{WeiXu}, which builds on \cite{CorrespondenceICMmatrices}. We have that $\GL_r(A_\p)$ acts on $\{M \in \Mat_r(A_\p) \mid \charpoly(M) = f\}$. We split the set into its $\GL_r(A_\p)$ orbits and then use the orbit stabilizer theorem. To count the number of orbits we use \th\ref{LFWeiXu} that relates the number of orbits to $\ICM(A_\p[x]/f(x))$. To count points in the stabilizers we first rewrite them in terms of orders in $K_\p = \Frac(A_\p[x]/f(x))$, which is done in Section \ref{secStabilizers}. In Section \ref{secLocalFactors} we do the counting and combine everything to get our result. 

\subsection{The Characteristic 0 Case}
In Section \ref{secChar0} we compare the function field setting to the number field setting. Regarding the computations of the ideal class monoid everything caries over. Gekeler ratios in characteristic 0 become different, so this does not carry over well. 

\subsection{Implementation} 
We have an implementation of the algorithm in Magma, see \cite{codeLocalGekelerRatios}. Some remarks regarding the code are collected in Section \ref{secImplementation}.

	\section{Acknowledgments}
	I thank Valentijn Karemaker for her supervision during this project. I thank Amie Bray for the calls we made about Gekeler ratios. 
	
	\section{Notation}
	For convenience we added a reference list of the used notation.\\
	We use the convention that overorders of $R$, $R_{(\p)}$ and $R_\p$ are always $A$-orders, $A_{(\p)}$-orders and $A_\p$-orders, respectively.
	\begin{itemize}
		\item $A = \F_q[T]$;
		\item $F = \F_q(T)$;
		\item $f(x) \in A[x]$ is a polynomial of degree $r$;
		\item $K = \Frac(A[x]/f(x))$ the fraction field of $A[x]$ modulo $f$;
		\item $R$ is an $A$-order in $K$;
		\item $\p \subset A$ is a prime ideal;
		\item $\p_1,\ldots,\p_k$ are the primes of $R$ above $\p$;
		\item $R_{\p_i}$, $R_\p$ and $A_\p$ are the completions of $R$ and $A$ at $\p_i$ and $\p$ respectively;
		\item $R_{(\p_i)}$ and $A_{(\p)}$ are the localisations of $R$ and $A$ at $\p_i$ and $\p$ respectively;
		\item $R_{(\p)} = A_{(\p)}+\ldots+a_{r-1}A_{(\p)}$, where $1,a_1,\ldots,a_{r-1}$ is an $A$-basis of $R$;
		\item $\ord(M)$ for a module $M$ is the order ideal, see \th\ref{OOdefiOrdIdeal};
		\item $[M:N]$ for modules $M$ and $N$ is the index ideal, see \th\ref{OOdefiOrdIdeal};
		\item $(I:J)$ is the ideal quotient or colon ideal of the fractional ideals $I$ and $J$;
		\item $I \sim J$ for two fractional ideals if they are weakly equivalent;
		\item $W(R)$ is the set of weak equivalence classes of fractional $R$-ideals;
		\item $[I]$ for a fractional ideal $I$ is the weak equivalence class of $I$;
		\item $\{I\}$ for a fractional ideal $I$ is the class of $I$ in the ideal class monoid;
		\item $[M]$ for a matrix $M$ is its orbit under conjugation;
		\item $[a]_I$ for $a\in A$ and $I\subset A$ is the equivalence class $a + I$ in $A/I$;
		\item $W_S(R) = \{[I] \in W(R)\mid (I:I)=S\}$;
		\item $\Mat_r(S)$ is the set of $r\ti r$ matrices over $S$;
		\item $M_0 \in \Mat_r(A)$ is the matrix in rational canonical form with characteristic polynomial $f$ as above, see \th\ref{CSdefM0};
		\item $\charpoly(M)$ for a matrix $M$ is the characteristic polynomial of $M$;
		\item $\IntCl(B,L)$ for rings $B$ and $L$ is the integral closure of $B$ in $L$;
		\item $v_\p(f) = \lim_{n \ra \infty} \frac{|\{M \in \Mat_r(A/\p^n)\mid \charpoly(M)=f\}}{|\SL_r(A/\p^n)|/|\p|^{n(r-1)}}$ is the Gekeler Ratio of $f$ at $\p$;
		\item $m_\p$ is the number of orbits in the set $\{M \in \Mat_r(A_\p)\mid \charpoly(M)=f\}$ under conjugation with $\GL_r(A_\p)$;
	\end{itemize}

\section{Ideal Classes} \label{secICM} 

In this section we will show how to compute overorders and ideal classes of orders of a completed ring. For this we first collect some lemmas we need later on. Next we follow the approach of \cite{Stefano}: we show how to compute overorders in Paragraph \ref{parICOO}, how to compute weak equivalence classes in Paragraph \ref{parICWE} and the ideal classes themselves in Paragraph \ref{parICICM}. 

\subsection{Some Results on Ring Structures and Picard Groups} \label{parIClemmas} \phantom{=}\\
In this paragraph we will collect some lemmas that we will need in the rest of this section. 

\begin{notation}
We recall the following notation: $f \in A[x]$, $K$ is the fraction field of $A[x]/f$ and $R$ is an $A$-order in $K$. Let $\p \subset A$ be a prime ideal and let $\p_1,\ldots,\p_k$ be the primes of $R$ above $\p$. We will denote the completions with $R_{\p_i}$, $R_{\p}$ and $A_{\p}$. We will denote the localisations with $R_{(p_i)}$ and $A_{(p)}$.
\end{notation} 

\begin{defi}
Let $1,a_1,\ldots,a_{r-1}$ be an $A$-basis of $R$. We set $R_{(\p)} = A_{(\p)} + a_1 A_{(\p)} + \ldots + a_{r-1} A_{(\p)}$.
\end{defi}

\begin{lemma}\th\label{spRpproduct} 
We have $R_\p \cong R_{\p_1} \ti \cdots \ti R_{\p_k}$.
\end{lemma} 

\begin{proof}
We have that $\p \subset \p_i$, which gives, since the $\p_i$ are coprime, that $\p \subset \cap_{i=1}^k \p_i = \prod_{i=1}^k \p_i$. From the Kummer-Dedekind theorem, \cite[Theorem 3.1]{dictaatANT}, we get that $\prod_{i=1}^k \p_i^{e_i} \subset \p$, where $e_i$ is the ramification index of $\p_i$. Together this gives us that for all $n\geq 1$ we have $\prod_{i=1}^k \p_i^{ne_i} \subset \p^n \subset \prod_{i=1}^k \p_i^n$, so we get surjective maps 
$$\prod_{i=1}^k R/\p_i^{ne_i} \ra R/\p^n \ra \prod_{i=1}^k R/\p_i^n .$$ 
Since the rings $\prod R/\p^{e_in}$ and $R/\p^n$ are finite for all $n$, we get by \cite{surjsystem} surjective maps when taking the inverse limit:
$$\varprojlim_n \left( \prod_{i=1}^k R/\p_i^{e_in}\right) \xrightarrow{\phi} \varprojlim_n R/\p^n \xrightarrow{\psi} \varprojlim_n \left( \prod_{i=1}^k R/\p_i^n \right).$$
Note that these inverse limits simplify to $\prod_{i=1}^k R_{\p_i} \xrightarrow{\phi} R_\p \xrightarrow{\psi} \prod_{i=1}^k R_{\p_i}$, since inverse limits and direct products commute. We will prove that $\psi \circ \phi$ is the identity map. We have
\begin{align*} 
	\psi(\phi(\prod_{i=1}^k ([a_{i,n}]_{\p_i^n})_{n\geq 1} )) 
	&= \psi(\phi(\prod_{i=1}^k ([a_{i,e_i n}]_{\p_i^{e_in}})_{n\geq 1}  )) \\
	&= \psi( ([(a_{1,e_1n},\ldots,a_{k,e_kn})]_{\p^n})_{n\geq 1} ) \\
	&= ([(a_{1,e_1n},\ldots,a_{k,e_kn})]_{\prod \p_i^n} )_{n \geq 1} \\
	&= ([a_{1,e_1n}]_{\p_1^n},\ldots,[a_{k,e_kn}]_{\p_k^n} )_{n \geq 1} \\   
	&= \prod_{i=1}^k ([a_{in}]_{\p_i^n})_{n\geq 1}. 
\end{align*} 
For the first equality we use that the congruence classes at the levels divisible by $e_i$ define the same element of the inverse limit as when we are using all the levels. For the second equality we use that the tuple $([a_{1,e_1n},\ldots,a_{k,e_kn}])$ defines with the Chinese remainder theorem a unique class in $R/\p^n$, which we denote with $[(a_{1,e_1n},\ldots,a_{k,e_kn})]_{\p^n}$. For the third equality we project onto $\prod_{i=1}^k \p_i^n$. In the fourth equality we project further and ``undo" our Chinese remainder theorem. For all $i$ it holds per construction that $[a_{1,e_1n},\ldots,a_{k,e_kn}]_{\p^n} \equiv a_{e_in} \mod \p_i^{e_in}$, so this equality also holds modulo $\p_i^n$. In the last equality we use that direct products and inverse limits commute.  
This proves $\psi \circ \phi = \text{id}$. Since we already know that $\phi$ is surjective, this gives that $\phi$ is a bijection with inverse $\psi$ and $\prod_{i=1}^k R_{\p_i} = R_\p $. 
\end{proof}

\begin{lemma}\th\label{Voight953} \emph{(\cite[9.5.3]{Voight})} Let $D$ be a discrete valuation ring with maximal ideal $\p$ and let $V$ be a finite dimensional vector space over the field of fractions of $D$. Then there is a bijection between the set of $D$-lattices $M$ in $V$ and the set of $D_\p$-latices $M_\p$ in $V_\p$ given by:
$$M \ra M_\p = M \te_{R_{(\p)}} R_\p \text{ and } M_\p \ra M_\p \cap V.$$
\end{lemma}

\begin{lemma} \th\label{spRplocalisation} 
We have $R_{(\p)} \cong R_{R\bs [\p_1\cup \ldots \cup \p_k]}$, the localisation of $R$ in the set $R\bs [\p_1 \cup \ldots \cup \p_k]$.
\end{lemma}

\begin{proof}
We have by \th\ref{Voight953} that $R_\p \cap K = R_{(\p)}$. Let $\psi: R_\p \ra R_{\p_1} \ti \cdot \ti R_{\p_k}$ be the bijection of \th\ref{spRpproduct}. \\
We will prove that $\psi(R) = \{ (r,r,\ldots,r) \in \prod_{i=1}^k R_{\p_i} \mid r \in R\}$. Let $s \in R$. We choose to work with its representative $([s]_{\p^n})_{n\geq 1}$. Then $\psi(s) = ( ([s]_{\p_1^n})_{n\geq 1},\ldots,([s]_{\p_k^n})_{n\geq 1}) = (s,\ldots,s)$. Since $\psi$ is well-defined, this gives that $\psi(R) = \{(r,\ldots,r)\in R_{\p_1}\ti \cdots \ti R_{\p_k} \mid r \in R\}$. \\
We have that $K = \Frac(R)$, so 
$$R_{(p)} = R_\p \cap K = (R_{\p_1}\ti \ldots \ti R_{\p_k}) \cap \{(r,\ldots,r) \mid r \in K\} \subseteq R_{(\p_1)} \ti \cdots R_{(\p_k)}.$$ 
We have that $R_{(\p_i)} = \{ \frac{a}{b} \mid a \in R,\ b \in R\bs \p_i\}$. This implies that $\frac{r}{s} \in R_\p \cap K$ if and only if it is equivalent to $\frac{a}{b}$ with $a \in R$, $b \in R \bs [\p_1 \cup \ldots \cup \p_k]$. The set $R\bs [\p_1 \cup \ldots \cup \p_k]$ is multiplicatively closed, so $R_{(\p)}$ is equal to the localisation of $R$ in this set. 
\end{proof}

\begin{lemma} \th\label{spPicTrivRpi} 
The Picard group of $R_{\p_i}$ is trivial.
\end{lemma}

\begin{proof}
Let $I \subset R_{\p_i}$ be an invertible ideal. We are going to prove that $I$ is principal. \\
Since $I$ is invertible, there exists a fractional $R_{\p_i}$ ideal $J$ such that $IJ = R_{\p_i}$. This means that there exist $x_j \in I$ and $y_j \in J$ such that $\sum_{j=1}^m x_jy_j = 1$. Note that all products $x_jy_j \in R_{\p_i}$, but they cannot be all in $\p_i$ since $1 \notin \p_i$. Without loss of generality $x_1y_1 \in R_{\p_i} \bs \p_i$. Since $\p_i$ is the only maximal ideal in $R_{\p_i}$, $x_1y_1$ is invertible. Now take $x\in I$. Then
$$x = x_1 \cdot x y_1 \cdot (x_1y_1)^{-1} \in x_1 R_{\p_i} \cdot R_{\p_i} = x_1 R_{\p_i}.$$
This implies that $I = x_1R_{\p_i}$, so it is principal. 
\end{proof}

\begin{cor} \th\label{spPicTrivRp} 
The Picard group of $R_\p$ is trivial.
\end{cor}

\begin{proof}
Since ideals and principal ideals split over products, we get with \th\ref{spRpproduct} that 
$$\Pic(R_{\p}) = \Pic(R_{\p_1} \ti \cdot \ti R_{\p_k}) = \Pic(R_{\p_1}) \ti \cdots \ti \Pic(R_{\p_k}) = 1\ti \cdots \ti 1 = 1.$$
\end{proof}

\subsection{Computing the Overorders} \label{parICOO}\phantom{=}\\
In this paragraph we will introduce $\p$-overorders and show that the overorders of $R_\p$ correspond with the $\p$-overorders of $R$ if $R= A[x]/f(x)$ for some irreducible $f$. 

\begin{defi} \th\label{OOdefiOrdIdeal} 
Let $D$ be a Dedekind domain and let $M$ be a finitely generated torsion module over $D$. Then by \cite[Theorem 6.3.20]{IntroRM} we get that $M \cong D/\a_1 \oplus \ldots \oplus D/\a_n$ for some $\a_1,\ldots,\a_n \in D$. The \textbf{order ideal} of $M$ is $\ord(M) = \prod_{i=1}^n \a_i$. If $N$ and $N'$ are both $D$-modules, $N' \subset N$ and $N/N'$ is a torsion module, then the \textbf{index ideal} is $[N:N'] = \ord(N/N')$.  
\end{defi}

\begin{defi}(\cite[Definition 3.2]{OverordersHS}) 
Let $R \subset S$ be a finite extension of $A$-orders. Then $S$ is a \textbf{$\p$-overorder} of $R$ if the index ideal $[R:S]$ is a $\p$-power. 
\end{defi}

\begin{prop} \th\label{OObijApandAorders} 
Assume that $f \in A[x]$ is irreducible and that $R = A[x]/f(x)$. Then the map 
$$\phi: \{\p\text{-overorders of }R\} \ra \{A_{(\p)}\text{-overorders of }R_{(\p)} \},\ S \ra S_{(\p)}$$
is a bijection. 
\end{prop}

\begin{proof}
We first prove surjectivity. Let $S'$ be an $A_{(\p)}$-overorder of $R_{(\p)}$. Let $\O$ be the maximal $\p$-overorder of $R$. Define $S = \{x \in \O \mid [x]_R \in S'/R_{(\p)} \}$. This is well-defined, since $\O/R \cong \O_{(\p))} / R_{(\p)} \supseteq S'/R_{(\p)}$. Since $\O$ is a $\p$-overorder of $R$, the same holds for $S$. We will prove that $\phi(S) = S'$. We get
\begin{align*}
	\phi(S) &= \{xy \mid x \in \O,\ y \in A_{(\p)},\ [x]_R \in S'/R_{(\p)} \} \\
	&= \{xy \mid x \in \O,\ y\in A_{(\p)},\ [xy]_{R_{(\p)}} \in S'/R_{(\p)} \} \\
	&= \{z \in \O_{(\p)} \mid [z] \in S'/R_{(\p)} \}. 
\end{align*} 
We have that $S'\subset \O_{(\p)}$, since $\O_{(\p)}$ is the maximal $A_{(\p)}$-order in $K$. This gives $S'\subset \phi(S)$. On the other hand $\phi(S) \subset S'+ R_{(\p)} = S'$, which gives $\phi(S) = S'$. \\
Now we move to injectivity. Let $S$ and $S'$ be $\p$-overrings of $R$ such that $\phi(S) = \phi(S')$. We have that $S/R \cong \phi(S)/R_{(\p)} = \phi(S')/R_{(\p)} \cong S'/R$ as $A$-modules. This gives that
$$S/R \cong A/\p^{\d_1} \oplus \ldots \oplus A/\p^{\d_{r-1}} \cong S'/R,$$
with $\d_1 \leq \d_2\leq \ldots \leq \d_{r-1}$. We only have $\p$-th powers, since $S$ and $S'$ are $\p$-overrings. \\
With \cite[Theorem 2.1]{EndDM} we get 
\begin{align*}
	S = A + \frac{f_1(\pi)}{\p^{\d_1}}A + \ldots + \frac{f_{r_1}(\pi)}{\p^{\d_{r-1}}} A, \\
	S' = A + \frac{g_1(\pi)}{\p^{\d_1}}A + \ldots + \frac{g_{r_1}(\pi)}{\p^{\d_{r-1}}} A, \\
\end{align*}
where $\pi$ is a zero of $f$ and $f_i,g_i$ are monic polynomials of degree $i$. We have with Gauss' lemma that $f$ is irreducible over $F$, so it is also irreducible over $A_{(\p)}$. This means that with the same theorem and \cite[Remark 2.2]{EndDM} we get that 
\begin{align*}
	S_{(\p)}  = A_{(\p)} + \frac{f_1(\pi)}{\p^{\d_1}}A_{(\p)} + \ldots + \frac{f_{r_1}(\pi)}{\p^{\d_{r-1}}} A_{(\p)}, \\
	S'_{(\p)} = A_{(\p)} + \frac{g_1(\pi)}{\p^{\d_1}}A_{(\p)} + \ldots +\frac{g_{r_1}(\pi)}{\p^{\d_{r-1}}} A_{(\p)}. \\
\end{align*}
With the same remark we get that if $\frac{g_i(\pi)}{\p^{\d_i}} \in S$ for all $i$, then $S = S'$. So suppose this is not the case and let $j$ be the smallest such that $\frac{g_j(\pi)}{\p^{\d_j}} \notin S$. We know that $S_{(\p)} = S'_{(\p)}$, which implies that $\frac{g_j(\pi)}{\p^{d_j}} \in S_{(\p)}$. This gives that there exist $a_0,\ldots,a_{r-1} \in A_{(\p)}$ such that 
\begin{align} \label{OOeqbij} 
	a_0 + \frac{f_1(\pi)}{\p^{\d_1}} a_1 + \ldots + \frac{f_{r-1}(\pi)}{\p^{\d_{r-1}}} a_{r-1} = \frac{g_j(\pi)}{\p^{\d_j}} 
\end{align}
Looking at the degree in $\pi$ gives that $a_{j+1}=\ldots = a_{r-1} = 0 \in A$. We will prove with induction on $i$ from $j$ to $0$ that $a_i \in A$ for all $i$.\\ 
\textit{Induction basis:} Comparing the coefficient of $\pi^j$ in \eqref{OOeqbij} gives since $g_j(\pi)$ and $f_j(\pi)$ are monic that $a_j = 1 \in A$. \\
\textit{Induction hypothesis:} Suppose that $a_k,a_{k+1},\ldots,a_j \in A$.\\
\textit{Induction step:} Let $f_{i,m}$ be the coefficient of $x^m$ in $f_i$ and the same for the $g_i$. If we compare the coefficient of $\pi^{k-1}$ in \eqref{OOeqbij} we get
$$\frac{a_{k-1}}{\p^{\d_{k-1}}} + \frac{a_k f_{k,k-1}}{\p^{\d_k}} + \ldots + \frac{a_jf_{j,k-1}}{\p^{\d_j}} = \frac{g_{j,k-1}}{\p^{\d_j}}.$$
Since with the induction hypothesis $f_{i,k-1},g_{j,k-1}$ and $a_m$ for $k \leq m \leq j$ are all in $A$ and $\d_{k-1} \leq \ldots \leq \d_j$ we get from this that $a_{k-1} = \frac{c}{\p^\d}$ with $c \in A$ and $\d \geq 0$ and without loss of generality we may assume that $\p \nmid c$. We know that $a_{k-1} \in A_{(\p)}$, so $0 \leq \nu_{\p}(a_{k-1}) = \nu_{\p}(c) - \nu_{\p}(\p^{\d}) = 0 - \d$. This implies that $\d = 0$ and $a_{k-1} \in A$ as needed. \\
This completes the induction, so we get that $\frac{g_j(\pi)}{\p^{\d_j}} \in S$. This is a contradiction with our assumption, so $S = S'$ and $\phi$ is injective.  
\end{proof}

\begin{lemma} \th\label{OObijApandAporders} 
There is a bijection between the $A_{(p)}$-overorders of $R_{(\p)}$ and the $A_\p$-overorders of $R_\p$ given by $S_{(\p)} \ra S_{(\p)}\te A_\p$ and $S_\p \ra S_\p \cap K$. 
\end{lemma}

\begin{proof}
We have that $A_{(\p)}$-overorders of $R_{(\p)}$ are $A_{(\p)}$-lattices in $K$ that are also rings, and $A_\p$-overorders of $R_\p$ are $A_\p$-lattices in $K_\p$. Restricting the bijection in \th\ref{Voight953} gives the result. 
\end{proof}

\subsection{Weak Equivalence Classes} \label{parICWE}\phantom{=}\\
In this paragraph we will show how to extend the definition of weak equivalence to the local setting and compute weak equivalence classes. 

\begin{prop}\emph{(c.f. \cite[Proposition 4.1]{Stefano})} \th\label{WEWeakEquiv} 
Let $I,J$ be fractional $R_\p$-ideals. Then the following are equivalent:
\begin{enumerate}
	\item For all $i \in \{1,\ldots,k\}$ there exists an $x \in K_{\p_i}$ such that $I_{\p_i} = xJ_{\p_i}$;
	\item $1 \in (I:J)(J:I)$;
	\item $(I:I) = (J:J) = S$ and there exists an invertible fractional $S$-ideal $L$ such that $I = LJ$. 
\end{enumerate}
\end{prop} 

\begin{proof}
The proof of $(1)\Ra (2)$ and $(2)\Ra (3)$ is the same as that of \cite[Proposition 4.1]{Stefano}. For $(3) \Ra (1)$, let $I,J$ be fractional $R_\p$-ideals with multiplicator ring $S$ and let $L$ be an invertible fractional $S$-ideal such that $I = LJ$. With \th\ref{spPicTrivRp} we have that the Picard group of $S$ is trivial, so $L = xS$ for some $x \in K$. This gives $I_{\p_i} = (xSJ)_{\p_i} = x (SJ)_{\p_i} = x ((J:J)J)_{\p_i} = xJ_{\p_i}$.   
\end{proof}

\begin{defi} \th\label{WEDefWeakEquiv} 
We say that two fractional $R_\p$-ideals $I$ and $J$ are \textbf{weakly equivalent} if the conditions of \th\ref{WEWeakEquiv} hold. We write $I \sim J$ for this. The weak equivalence class of $I$ is denoted $[I]$. We denote the set of weak equivalence classes by $W(R)$ and we denote $W_{S}(R) = \{[I] \in W(R) \mid (I:I) = S\}$.   
\end{defi} 

In the following series of lemmas we will show that we can decide whether or not $I_\p \sim J_\p$ by computations in $R$.  

\begin{lemma}\th\label{WEequivLocCompl} 
Let $I,J$ be fractional ideals $R$-ideals and let $\p \subset A$ be prime. Then $I_\p \sim J_\p$ if and only if $I_{(\p)} \sim J_{(\p)}$. 
\end{lemma} 

\begin{proof}
Suppose that $I_{(\p)} \sim J_{(\p)}$. Then there exists an $S$ such that $(I_{(\p)}:I_{(\p)}) = (J_{(\p)}:J_{(\p)}) = S$ and there exists an invertible $S$-ideal $L$ such that $I_{(\p)} = LJ_{(\p)}$. We have $(L_{(\p)}\cdot J_{(\p)})_\p = L_\p\cdot J_\p = I_\p$ and $(I_{(\p)}:I_{(\p)})_\p = S_\p$ and $(L\cdot L^{-1})_\p = L_p \cdot L^{-1}_\p = S_\p$. This implies that $L_\p$ is an invertible $S_\p$-ideal such that $I_\p = L_\p \cdot J_\p$, so $I_\p \sim J_\p$. \\
Suppose now that $I_\p \sim J_\p$. Then $(I_\p:I_\p) = (J_\p:J_\p) = S_\p$ and there exists some invertible $S_\p$ ideal $L_\p$ with inverse $L'_\p$ such that $I_\p = L_\p J_\p$. By \th\ref{Voight953} we have that the ideals $I_{(\p)}$, $J_{(\p)}$, the ring $S_{(\p)}$ and the $S_{(\p)}$-ideals $L_{(\p)}$ and $L'_{(\p)}$ exist. First $(I_{(\p)}:I_{(\p)})_\p = (I_\p:I_\p) = S_\p = (S_{(\p)})_\p$ implies $(I_{(\p)}:I_{(\p)}) = S_{(\p)}$ and similar $(J_{(\p)}:J_{(\p)}) = S_{(\p)}$. Then $(L_{(\p)}\cdot J_{(\p)})_\p = L_\p\cdot J_\p = I_\p = (I_{(\p)})_\p$ and $(L_{(\p)}\cdot L'_{(\p)})_\p = L_p \cdot L'_\p = S_\p = (S_{(\p)})_\p$ imply that $L_{(\p)}$ is an invertible $S_{(\p)}$ ideal such that $L_{(\p)}J_{(\p)} = I_{(\p)}$. This gives $I_{(\p)} \sim J_{(\p)}$.
\end{proof}

\begin{lemma}\th\label{WEcalcLocWE} 
Let $I$ and $J$ be fractional $R$-ideals. Then $I_{(\p)} \sim J_{(\p)}$ if and only if 
$$((I:J)\cdot (J:I)) \cap (R \bs [\p_1\cup \ldots \cup \p_k]) \not= \emptyset.$$
\end{lemma}

\begin{proof}
With \th\ref{WEWeakEquiv} we have that $I_{(\p)} \sim J_{(\p)}$ if and only if $1 \in (I_{(\p)}:J_{(\p)}) \cdot (J_{(\p)}:I_{(\p)})$. \\
In the same way as \th\ref{spRplocalisation} we get that $I_{(\p)}$ and $J_{(\p)}$ are the localisations of $I$ and $J$ in $R \bs [\p_1\cup \ldots \p_k]$. This gives that
$$(I_{(\p)}:J_{(\p)})\cdot (J_{(\p)}:I_{(\p)}) = (I:J)_{(\p)} \cdot (J:I)_{(\p)} = [(I:J)\cdot (J:I)]_{(\p)}.$$
If $(I:J)(J:I) \cap R\bs[\p_1\cup \ldots \cup \p_k] \not= \emptyset$, say it contains $s$, then $\frac{s}{s} = 1 \in [(I:J)(J:I)]_{(\p)}$. On the other hand, if $1 \in [(I:J)(J:I)]_{(\p)}$, then $1 = \frac{s}{s}$ with $s \in (I:J)(J:I)$ and $s \in R\bs [\p_1 \cup \ldots \cup \p_k]$. 
\end{proof}

Next we show that we can compute $W_{S_\p}(R_\p)$ from $W_S(R)$. With \cite[Algorithm 3]{Stefano} we can compute $W_S(R)$.

\begin{lemma} \th\label{WEmapidealsurj} 
The map $\chi: \{\text{fractional ideals in }S\} \ra \{\text{fractional ideals in }S_\p\}$, $I \ra I_\p$ is surjective.
\end{lemma}  

\begin{proof}
We can restrict the bijection of \th\ref{Voight953} to ideals, since non-zero ideals are lattices and the bijection preserves being an ideal. This gives a bijection between the set of ideals in $S_{(\p)}$ and the set of ideals in $S_\p$ given by $I_{(\p)} \ra I_\p$. With \th\ref{spRplocalisation} we get that $S_{(\p)} = S_{S\bs[\p_1\cup \ldots \cup \p_k]}$, so with \cite[Lemma 2.8]{dictaatANT} every ideal $I_{(\p)}$ in $S_{(\p)}$ is of the form $(S \bs [\p_1\cup \ldots \cup\p_k])^{-1} I'$ for some ideal $I' \subset S$. We have that $\chi(I')$ is the ideal corresponding to $I_{(\p)}$, i.e. $I_\p = \chi(I')$.  
\end{proof}

\begin{lemma} \th\label{WEcompsurjonWE} 
The map $\theta: W(R) \ra W(R_\p)$ defined by $[I] \ra [I_\p]$ is well-defined and surjective. 
\end{lemma}

\begin{proof}
Suppose that $I \sim J$. Then $1 \in (I:J)(J:I)$ by \th\ref{WEWeakEquiv}. This gives 
$$1 \in (I:J)(J:I) \subset [(I:J)(J:I)]_{(\p)} = (I_{(\p)}:J_{(\p)})(J_{(\p)}:I_{(\p)}).$$
With \th\ref{WEcalcLocWE} this gives that $I_{(\p)} \sim J_{(\p)}$. Then \th\ref{WEequivLocCompl} gives $I_\p \sim J_\p$, so $\theta$ is well-defined.\\
Let $[J] \in W(R_\p)$ and let $S_\p = (J:J)$. Then by \th\ref{WEmapidealsurj} there exists an overorder $S$ of $R$ and an $S$-ideal $I$ such that $I_\p = J$. Then $S_\p$ is an $S$-overorder, so also an $R$-overorder. This gives that $[I] \in W(R)$ and $\theta([I]) = [J]$. 
\end{proof}

\begin{lemma} \th\label{WEcompsurjonWES} 
Suppose that $f$ is irreducible and that $R = A[x]/f(x)$. Let $S$ be the $\p$-overorder of $R$ corresponding to $S_\p$ under the map of \th\ref{OObijApandAorders}. Then the map $\theta: W_S(R) \ra W_{S_\p}(R_\p)$ defined by $[I] \ra [I_\p]$ is well-defined and surjective.
\end{lemma}

\begin{proof}
Let $[I] \in W_S(R)$. Then $(I:I) = S$, so $(I_\p:I_\p) = (I:I)_\p = S_\p$, which gives that $[I_\p] \in W_{S_\p}(R_\p)$.  This implies that we can restrict the map of \th\ref{WEcompsurjonWE}.\\
Take $[\tilde{J}] \in W_{S_\p}(R_\p)$. Then by \th\ref{Voight953} there is some $J_{(\p)} \subset S_{(\p)}$ such that its completion is $\tilde{J}$. With the local-global dictionary for lattices, \cite[9.4.9]{Voight}, there exists a lattice $I$ such that $I_{(\p)} = J_{(\p)}$ and $I_{(\q)} = S_{(\q)}$ for all other primes $\q$ of $A$. Since $J_{(\p)}$ and $S_{(\q)}$ are ideals we get that $I$ is also an ideal. We have $(I:I)_{(\q)} = (I_{(q)}:I_{(\q)}) = (S_{(\q)}:S_{(\q)}) = S_{(q)}$ and $(I:I)_{(\p)} = (J_{(\p)}:J_{(\p)}) = S_{(\p)}$. Together this implies, since an order is determined by its localisations, that $(I:I) = S$, so $[I] \in W_S(R)$. By construction $\psi([I]) = [\tilde{J}]$. 
\end{proof}

\begin{remark} \th\label{ICMWeakEquivInsep}
In the case that $f$ is inseparable, we would like to calculate the weak equivalence classes with the same algorithms as in \cite{Stefano}. Most of the algorithm still works, but there is one lemma that does not work. This is because if $f$ is inseparable, the trace form is degenerate. This means that the statement of \cite[Proposition 5.1]{Stefano} does not make sense. We do have a variant of it with a similar proof, see below. We can substitute the lemma below for \cite[Proposition 5.1]{Stefano} and then we get a working algorithm in the inseparable case. 
\end{remark}

\begin{lemma}
Let $S$ be an overorder of $R$. Then every class in $W_S(R)$ has a representative $I$ satisfying $(S:\O_K) \subset I \subset \O_K$. 
\end{lemma}

\begin{proof}
Let $[I']$ be a class in $W_S(R)$. This means that $(I':I') = S$. Every ideal in $\O_K$ is invertible, so $I'\O_K$ is invertible. Its inverse $(\O_K:I'\O_K)$ is also invertible. Define $\phi: \Pic(S) \ra \Pic(\O_K)$ by $L \ra L\O_K$. This map is surjective, see \cite[Corollary 2.1.11]{DadeTausskyZassehausTheoryOfOrders}. Let $J$ be such that $\phi(J) = (\O_K:I'\O_K)$ and set $I = I'J$. We have that $[I'] = [I]$, since $J$ is invertible in $S$. Furthermore $I\O_K = I'J\O_K = I'(\O_K:I'\O_K) = \O_K$, so $I \subset \O_K$. We also have $(S:\O_K) \O_K = (S:\O_K)$, so $(S:\O_K) I = (S:\O_K) \O_K I = (S:\O_K) \O_K = (S:\O_K)$. Since $(S:\O_K) \subset S = (I:I)$ we get $(S:\O_K) = (S:\O_K)I \subset (I:I)I \subset I$. This implies that $I$ is the representative we need.  
\end{proof}

\subsection{Computing the Ideal Class Monoid} \label{parICICM} 

\begin{defi}
The ideal class monoid of $R$ is the set of fractional ideals of $R$ modulo the principal ideals. We write $\ICM_S(R) = \{\ \{I\} \in \ICM(R) \mid (I:I) = S\}$ for the subset of the ideal class monoid with ideals with multiplicator ring $S$. 
\end{defi} 

\begin{lemma} \th\label{spICMunion}
We have $\ICM(R_\p) = \sqcup_{\tilde{S}} \ICM_{\tilde{S}}(R_\p)$ where the union is over the $A_\p$-overorders of $R_\p$.  
\end{lemma}

\begin{prop}\emph{(c.f. \cite[Theorem 4.6]{Stefano})} \th\label{spICMisPicWE} 
Let $\tilde{S}$ be an overorder of $R_\p$. We have that
$$\ICM_{\tilde{S}}(R_\p) = \{ \{I_i\cdot J_j\} \mid I_i \in W_{\tilde{S}}(R),\ J_j \in \Pic(\tilde{S}) \}$$
and if $\{I_iJ_j\} = \{I_{i'}J_{j'}\}$ then $i=i'$ and $j=j'$. 
\end{prop}

\begin{proof}
We can copy the proof of \cite[Theorem 4.6]{Stefano}, where we replace \cite[Proposition 4.1]{Stefano} with \th\ref{WEWeakEquiv}. 
\end{proof}

\begin{prop} \th\label{ICMisWE} 
We have that $\ICM(R_\p) = \sqcup_{\tilde{S}} W_{\tilde{S}}(R_\p)$, where the union is taken over all $A_\p$-overorders of $R_\p$. 
\end{prop}

\begin{proof}
Let $\tilde{S}$ be an overorder of $R_\p$. From \th\ref{spPicTrivRp} we get that $\Pic(\tilde{S}) = \{\tilde{S}\}$. If $[\tilde{I}] \in W_{\tilde{S}}(R_\p)$, then $\tilde{S} = (\tilde{I}:\tilde{I})$ and we get that $\tilde{S}\tilde{I} = (\tilde{I}:\tilde{I})\tilde{I} = \tilde{I}$. With \th\ref{spICMisPicWE} we get the result.  
\end{proof}

\begin{thm} \th\label{ICMcomputable} 
Let $f \in A[x]$ be irreducible and let $R = A[x]/f$. We can compute $\ICM(R_\p)$. 
\end{thm} 

\begin{proof} \phantom{=}
\begin{itemize}
	\item We first compute the overorders $S$ of $R$ corresponding to the overorders $S_\p$ of $R_\p$. This can be done by \th\ref{OObijApandAorders,OObijApandAporders}. 
	\item Then we compute $W_S(R)$ for all overorders $S$ of $R$ of the previous step. This can be done with the algorithms in \cite{Stefano}, keeping in mind \th\ref{ICMWeakEquivInsep}.
	\item Then we can check with \th\ref{WEcalcLocWE,WEequivLocCompl} which classes coincide if we map from $W_S(R)$ to $W_{S_\p}(R_\p)$. This gives all classes with \th\ref{WEcompsurjonWES}. 
	\item Finally we take the disjoint union of the weak equivalence classes, since with \th\ref{ICMisWE} we get that $\ICM(R_\p) = \sqcup_{S_\p} W_{S_\p}(R_\p)$.  
\end{itemize}
\end{proof}

\begin{lemma} \th\label{rpICMtriv} 
Let $\p$ be a regular prime of $R$. Then $\ICM(R_\p)$ is trivial. 
\end{lemma}

\begin{proof}
We have with \th\ref{spRpproduct} that $R_\p \cong R_{\p_1} \ti \ldots R_{\p_k}$. Since $\p$ is regular, each $R_{\p_i}$ is a DVR. This implies that each ideal is principal, so $\ICM(R_{\p_i})$ is trivial. Since ideals and principal ideals split over products, we get that $\ICM(R_\p)$ is also trivial.  
\end{proof}

\section{Calculating Stabilizers} \label{secStabilizers}

In this section we will rewrite $\Stab_{\GL_r(A_\p)}(M_0)$ for a matrix $M_0$ in rational canonical form with characteristic polynomial $f$ to $\prod_{\q \mid \p} \O_\q$ where $\q$ runs over the primes of $K = \Frac(A[x]/f(x))$ above $\p$. We also have \th\ref{CSstabM1} for matrices not in rational canonical form. 

\begin{notation} \th\label{CSdefM0} 
Let $M_0 \in \Mat_r(A)$ be the $r\ti r$ matrix in rational canonical form with characteristic polynomial $f$. In other words, if $f = -b_1-b_2x-\ldots - b_{r}x^{r-1} + x^r$, then 
$$M_0 = \mat{0 & 0 & \ldots & 0 & b_1 \\ 1 & 0 & \ldots & 0 & b_2 \\ \vdots & \vdots & & \vdots & \vdots \\ 0 & 0 & \ldots & 1 & b_{r} \\} .$$  
\end{notation}

\begin{lemma} \th\label{CSfnotramified} 
If $f \in A[x]$ is irreducible over $A$ and $\p$ is prime, then $f \not= g^2 h$ for $g,h \in F_\p[x]$. 
\end{lemma}

\begin{proof}
Write $f = f_1\cdots f_n$ with $f_i \in F_{\p}[x]$ irreducible. The fact that $f$ is irreducible over $A$ means that $f$ defines a finite field extension of $F$, so we may use \cite[Theorem 2.8.5]{DM}. The proof of this theorem tells us that the $f_i$ are distinct, even if $f$ is not separable. This implies that $f \not= g^2h$ for $g,h \in F_\p[x]$.  
\end{proof} 

\begin{lemma} \th\label{CSMM0=M0M} 
Suppose that $f\not= g^2h$ for $g,h\in F_\p[x]$. We have $\{M \in \Mat_r(F_\p) \mid MM_0 = M_0 M\} = F_\p(M_0) \cong F_\p[x]/f(x)$. 
\end{lemma}

\begin{proof}
It is clear that $F_\p(M_0) \subset \{M \in \Mat_r(F_\p) \mid MM_0 = M_0M\}$. We note that $\dim_{F_\p}(F_\p(M_0)) = \deg(f) = r$, since $f$ is the largest invariant factor of $M_0$, so the minimal polynomial of $M_0$, \cite[Proposition 12.13]{DummitFoote}. Let 
$$M_0 = \mat{0 & 0 & \ldots & 0 & b_1 \\ 1 & 0 & \ldots & 0 & b_2 \\ \vdots & \vdots & & \vdots & \vdots \\ 0 & 0 & \ldots & 1 & b_r \\}  \text{ and } M = \mat{a_{1,1} & \ldots & a_{1,r} \\ a_{2,1} & \ldots & a_{2,r} \\ \vdots & & \vdots \\ a_{r,1} & \ldots & a_{r,r} \\}.$$
If we expand we get
$$M_0M = \mat{b_1a_{r,1} & \ldots & b_1a_{r,r} \\ a_{1,1}+b_2a_{r,1} & \ldots & a_{1,r}+b_2a_{r,r} \\ \vdots & &\vdots \\ a_{r-1,1} + b_ra_{r,1} & \ldots & a_{r-1,r}+b_ra_{r,r} }$$
and 
$$MM_0 = \mat{a_{1,2} & a_{1,3} & \ldots & a_{1,r} & \sum_{i=1}^r a_{1,i}b_i \\ \\\vdots &\vdots &  & \vdots \\ a_{r,2} & a_{r,3}& \ldots & a_{r,r} & \sum_{i=1}^r a_{r,i}b_i \\}.$$ 
The equality of these matrices is entrywise, which gives a system of $r^2$ linear equations in the $r^2$ variables $a_{i,j}$. We will prove that this system has rank $\geq r^2-r$. This will give us that the space of solutions has dimension at most $r$, so the set of solutions is $F_\p(M_0)$. \\
The subsystem of equations we will use is 
\begin{align*}
	&a_{1,i}-b_1a_{r,i} = 0 \ &\va\ 2\leq i \leq r \tag{S.i}\\
	&a_{i,j} - a_{i-1,j-1} - b_ia_{r,i-1} = 0 \ &\va\ 2 \leq i \leq r,\ 2 \leq j \leq r \tag{S.ij} 
\end{align*}
Order the variables lexicographically. We will prove with induction on $i$ that we can add for all $1 \leq i \leq r-1$ and for all $i < j \leq r$ an equation of the form $a_{i,j} + \sum_{k=1}^r c_{i,j,k} a_{r,k} = 0$ to the system without changing the set of solutions. \\
\textit{Induction basis:} For $i=1$ we have Equation $(S.j)$: $a_{1,j}-b_1a_{r,j} = 0$ already in the system for $1 < j \leq n$. \\
\textit{Induction hypothesis:} Suppose we can add equations of the form $a_{i,j} + \sum_{k=1}^r c_{i,j,k} a_{r,k} = 0$ to the system for a certain $i$ and all $i<j\leq r$. \\
\textit{Induction step:} For $(i+1,j)$ we have Equation $(S.(i+1)\ j)$ in the system:
$$-a_{i,j-1} + a_{i+1,j} - b_{i+1}a_{r,i} = 0$$
and with the induction hypothesis we have an equation
\begin{align} 
	a_{i,j-1}+\sum_{k=1}^r c_{i,(j-1),k} a_{r,k} \tag{IH.i\ (j-1)}.
\end{align} 
Adding these two together gives
$$a_{i+1,j} -b_{i+1}a_{r,i} + \sum_{k=1}^r c_{i,j,k}a_{r,k},$$ 
which is an equation of the required form that can be added to the system without changing the set of solutions. This proves the induction step.\\
Now we have for all $i<j$ an equation that starts with $a_{i,j}$ and has only bigger variables (in the lexicographical ordering) after that. \\
For $r-1 \geq i \geq j \geq 1$ we have Equation $(S,(i+1)\ (j+1))$ that starts with $-a_{i,j}$ and has only bigger variables after that. If we order these $r(r-1) = r^2-r$ equations by starting variable, we get an uppertriangular matrix with $1$'s and $-1$'s on the diagonal. So this subsystem has rank $r^2-r$, which proves that the rank of the complete system is at least $r^2-r$. 
\end{proof}

\begin{example}
Let $A = \F_3[T]$ and $f = 1-x+x^3$. This gives us $b_1=-1$, $b_2=1$ and $b_3=0$. We have 6 equations in our system, namely $(S.2)$ and $(S.3)$, which form the induction basis, and $(S.22),(S.23),(S.32)$ and $(S.33)$. We will do the induction step for $(i+1,j) = (2,3)$, which is the only case not covered by the induction basis. With the induction hypothesis for $(1,2)$ we get $a_{1,2} + \sum_{k=1}^3 c_{1,2,k}a_{3,k} = 0$. From $(S.2)$ we see that $c_{1,2,1}=c_{1,2,3} = 0$ and $c_{1,2,2}=1$ in this case. Equation $(S.23)$ gives us $a_{2,3}-a_{1,2}-a_{3,1} = 0$. Adding those two together gives $a_{2,3}-a_{3,1}+a_{3,2} = 0$. This is in the required form, with $c_{2,3,1} = -1$, $c_{2,3,2} = 1$ and $c_{2,3,3}=0$. 
\end{example}

\begin{lemma} \th\label{CSsetisintcl} 
Suppose that $f\not= g^2h$ for $g,h \in F_\p[x]$. We have that $\{M \in \Mat_r(A_\p) \mid MM_0 = M_0M\} = \IntCl(A_\p, F_\p(M_0))$, the integral closure of $A_\p$ in $F_\p(M_0)$.
\end{lemma}

\begin{proof}
First note that $\Mat_r(A_\p)$ is the maximal $A_\p$-order in $\Mat_r(F_\p)$. This implies that $\{M \in \Mat_r(A_\p) \mid MM_0 = M_0M\}$ is the maximal $A_\p$-order in $\{M \in \Mat_r(F_\p) \mid MM_0 = M_0M\}$. This last set is $F_\p(M_0)$ by \th\ref{CSMM0=M0M}. Since taking the integral closure of $A_\p$ in $F_\p(M_0)$ gives the maximal $A_\p$-order in $F_\p(M_0)$, we get $\{M \in \Mat_r(A_\p) \mid MM_0 = M_0M\} = \IntCl(A_\p,F_\p(M_0))$. 
\end{proof}

\begin{lemma} \th\label{CSintclisprodvalrings} 
Suppose that $f$ is irreducible over $A$ and let $K_\p = K \te F_\p$. Then $\IntCl(A_\p,K_\p) = \prod_{\q \mid \p} \O_\q$, where $\q$ are the primes of $K$ above $\p$ and $\O_\q$ is the valuation ring of the completion $K_\q$.
\end{lemma}

\begin{proof}
By \cite[Theorem 2.8.5]{DM} we get $K_\p = \prod_{w \mid \ord_\p} K_w$ where $w$ is a valuation extending $\ord_\p$. Since the primes above $\p$ correspond to the valuations extending $\ord_\p$, we get $K_\p = \prod_{\q \mid \p} K_\q$. The maximal $A_\p$-order in $K_\p$ is $\IntCl(A_\p,K_\p)$. The maximal $A_\p$-order in $\prod_{\q\mid\p}K_\q$ is $\prod_{q\mid\p}\O_\q$. Since maximal orders are unique, we get $\IntCl(A_\p,K_\p) \cong \prod_{\q\mid\p} \O_\q$. 
\end{proof}

\begin{lemma} \th\label{CSstabis} 
Let $f$ be irreducible over $A$. Let $\GL_r(A_\p)$ act on $\Mat_r(A_\p)$ by conjugation. We have $$\Stab_{\GL_r(A_\p)}(M_0) = \IntCl(A_\p,F_\p(M_0))^\ti = \left(\prod_{\q\mid\p} \O_\q\right)^\ti.$$  
\end{lemma}

\begin{proof}
We have that 
\begin{align*} 
\Stab_{\GL_r(A_\p)}(M_0) &= \{M \in \GL_r(A_\p) \mid MM_0 = M_0M\} \\
&= \{M \in \Mat_r(A_\p)\mid MM_0 = M_0M\}^\ti.
\end{align*}  
With \th\ref{CSfnotramified} we have that $f\not= g^2h$ for $g,h \in F_\p[x]$. With \th\ref{CSsetisintcl,CSintclisprodvalrings} we get the result. 
\end{proof}

\begin{lemma} \th\label{CSstabM1} 
Suppose that $f \not= g^2h$ for $g,h \in F_\p[x]$. Let $M_1 \in \Mat_r(A_\p)$ be such that its characteristic polynomial is $f$. We do not assume here that $M_1$ is in rational canonical form. Then there exists an $N \in \GL_r(F_\p)$ such that 
$$N^{-1}\IntCl(A_\p,F_\p(M_0))^\ti N = \Stab_{\GL_r(A_\p)}(M_1).$$ 
\end{lemma}

\begin{proof}
We have that $M_0$ and $M_1$ are in the same $\GL_r(F_\p)$-orbit, since they have the same characteristic polynomial, which determines a unique invariant factor. This means there exists some $N \in \GL_r(F_\p)$ such that $M_1 = NM_0N^{-1}$. This gives $\Stab_{\GL_r(F_\p)}(M_0) = N\Stab_{\GL_r(F_\p)}(M_1)N^{-1}$, so by \th\ref{CSMM0=M0M} that $N^{-1}F_\p(M_0)^\ti N = \Stab_{\GL_r(F_\p)}(M_1)$. \\
We will prove that $\IntCl(A_\p,N^{-1}F_\p(M_0)N) = N^{-1} \IntCl(A_\p,F_\p(M_0))N$.\\
Let $\a \in \IntCl(A_\p,F_\p(M_0))$. Then there exist $a_i \in A_\p$ such that $\a^r + a_{r-1}\a^{r-1}+\ldots + a_0 = 0$. Since $a_i \in A_\p \subset Z(\GL_r(F_\p))$, we get that $N^{-1}a_iN = a_iN^{-1}N = a_i$. This gives $(N^{-1}\a N)^r + a_{r-1}(N^{-1}\a N)^{r-1}+\ldots + a_0 = N^{-1} \a^r N + N^{-1} a_{r-1} \a^{r-1} N + \ldots + N^{-1} a_0 N = N^{-1}\cdot 0 \cdot N = 0$, which implies that $\a \in N^{-1}\IntCl(A_\p,F_\p(M_0))N$. In the same way we get the other inclusion. We have that $\Stab_{\GL_r(A_\p)}(M_1) = \IntCl(A_\p,F_\p(M_0))^\ti$, as $\Mat_r(A_\p)$ is the maximal order in $\Mat_r(F_\p)$. Together this gives $\Stab_{\GL_r(A_\p)}(M_0) = N^{-1}\IntCl(A_\p,F_\p(M_0))^\ti N$.   
\end{proof}


\section{Computing the Gekeler Ratios and their Product} \label{secLocalFactors} 
In this section we will show how a local Gekeler ratio can be computed and how to compute their product. We will first go from the ideal class monoid to orbits of the set $\{M \in \Mat_r(A_\p) \mid \charpoly(M) = f\}$. Then we use this to compute a single Gekeler ratio. Using that regular primes are easier than singular primes we provide an algorithm to compute the product of all (local) Gekeler ratios. 

\begin{lemma}\emph{(c.f. \cite[Lemma 6.4]{WeiXu})} \th\label{LFWeiXu} 
Let $S$ be a principal ideal domain and let $f(x) \in S[x]$ be a polynomial of degree $r$, such that there exist no $g,h \in S[x]$ such that $f=g^2h$. Let $\GL_r(S)$ act on $\Mat_r(S)$ by conjugation. Then there is a bijection between the number of orbits in the set $\{M \in \Mat_r(S)\mid \charpoly(M)=f\}$ and the number of ideal classes in $S[x]/f(x)$.  
\end{lemma} 

\begin{proof}
The difference between our formulation and the formulation of \cite[Lemma 6.4]{WeiXu} is that we have the condition $f\not= g^2h$, while they have the condition that $(f,f') = 1$. Denote with $L$ the fraction field of $S$.\\
The first time the condition $(f,f')$ is used in \cite[Lemma 6.4]{WeiXu} is to conclude that every matrix with characteristic polynomial $f$ is conjugated over $\GL_r(L)$ to the matrix $M_0$. This is true if there is one matrix in rational canonical form with characteristic polynomial $f$, which our condition also guarantees. \\
The other place in \cite[Lemma 6.4]{WeiXu} where the condition is used is to conclude that $\{M \in \Mat_r(L) \mid MM_0 = M_0M\} = L[x]/f(x)$. This is \th\ref{CSMM0=M0M}, since the proof works over any field.
\end{proof}

\begin{prop} \th\label{ICMsingprimes} 
We can compute the singular primes of $R = A[x]/f(x)$, even if $f$ is inseparable. 
\end{prop}

\begin{proof}
If $f$ is separable, the discriminant $D$ of $f$ is nonzero. All singular primes divide $D$, which gives a finite list to check. With the Kummer-Dedekind theorem, see \cite[Theorem 3.1]{dictaatANT}, we can then check which of these primes are singular. \\
If $f$ is inseparable, we get $D=0$, so this does not work. Instead we look at the primary decomposition of the conductor $\c = (R:\O_K)$. By \cite[Theorem 3.6]{ConradConductor} every prime ideal that is relatively prime to the conductor is invertible, so regular. This means that if $\p$ is singular, then $(\c,\p) \not = 1$. Since all nonzero prime ideals are maximal, this implies that $\c \subset \p$. Since $\c \not= 0$, we have that $\p$ is a minimal prime ideal containing $\c$. Those ideals correspond to the minimal elements in the set of radicals of primary ideals in the decomposition of $\c$.\\
Let $I \subset A[x]$ be an ideal containing $f$ and let $\bar{I}$ be its projection onto $A[x]/f(x)$. If $J_1,\ldots,J_n$ is a primary decomposition for $I$, then $\bar{J_1},\ldots,\bar{J_k}$ is a primary decomposition for $\bar{I}$. In $A[x] = \F_q[T,x]$ we can compute primary decompositions and their corresponding prime ideals with \cite{PrimaryDecompositionGTZ}. \\
This means we can compute the prime ideals that contain $\c$, which is a finite list. We then use Kummer-Dedekind to check which of these are indeed singular primes.  
\end{proof}

\begin{lemma} \th\label{LFcomporbits} 
Let $f(x) \in A[x]$ be irreducible. We can compute the number of orbits of $\{M \in \Mat_r(A_\p)\mid \charpoly(M)=f\}$ under the action of $\GL_r(A_\p)$ by conjugation for all primes $\p$ of $A$. 
\end{lemma}

\begin{proof}
We note that by \th\ref{CSfnotramified} $f$ is not of the form $g^2h$ for $g,h \in A_\p$. This means we can use \th\ref{LFWeiXu} with $S = A_\p$.\\
For all regular primes we have that $\ICM(R_\p)$ is trivial by \th\ref{rpICMtriv}. There are only finitely many singular primes of $R$ and we can compute which primes are singular with \th\ref{ICMsingprimes}. Since we can compute the number of ideal classes by \th\ref{ICMcomputable} for all singular primes, we also get the number of orbits.  
\end{proof}

\begin{defi}
The Gekeler ratio at the prime $\p$ is 
$$v_\p(f) = \lim_{n \ra \infty} \frac{|\{M \in \Mat_r(A/\p^n)\mid \charpoly(M)=f\}|\cdot |\p|^{n(r-1)}}{|\SL_r(A/\p^n)|}.$$
\end{defi}

\begin{remark}
Gekeler ratios are named after Gekeler, who described them for elliptic curves in \cite{GekelerEllipticCurves}. They come from a heuristic model that relates the Frobenius to a class in $\GL_2(\Z/n\Z)$. For Drinfeld modules (of rank $2$) Gekeler ratios are introduced in \cite[Section 7]{GekelerFrobeniusDistrDM}.
\end{remark}

Before we can compute these ratios, we need to lemmas that compute some limits.

\begin{lemma} \th\label{MGLandSL} 
We have $\lim_{n\ra \infty} \frac{|\GL_r(A/\p^n)|}{|\p|^{r^2n}} = (1-\frac{1}{|\p|}) \lim_{n\ra\infty} \frac{|\SL_r(A/\p^n)|}{|\p|^{(r^2-1)n}}$. 
\end{lemma}

\begin{proof}
We have an $|(A/\p^n)^\ti|$ to $1$ map $\phi_n: \GL_r(A/\p^n) \ra \SL_r(A/\p^n)$ defined by $M \ra \frac{1}{\Det(M)} M$. All elements in $A/\p^n$ that are not divisible by $\p$ are invertible, so $|(A/\p^n)^\ti| = |\p|^{n-1}(|\p|-1)$. This gives
\begin{align*}
	\lim_{n \ra \infty} \frac{|\GL_r(A/\p^n)|}{|\p|^{r^2n}}
	&= \lim_{n\ra\infty} \frac{|\p|^{n-1}(|\p|-1)|\SL_r(A/\p^n)|}{|\p|^{nr^2-n}|\p|^{n}} \\
	&= (1-\frac{1}{|\p|})\ \lim_{n\ra\infty} \frac{|\SL_r(A/\p^n)|}{|\p|^{n(r^2-1)}}, \\
\end{align*}
as required.
\end{proof} 

\begin{notation}
We write $N_K(\q) = |\O_K/\q|$ for the norm of $\q$ over $K$. For the norm over $F$ of an ideal $I \subset A$ we write $N_F(I)$ or $|I|$. 
\end{notation} 

\begin{lemma} \th\label{MprodOq} 
Let $\O_K$ be the maximal order in $K$, let $\q$ be a prime above $\p$ and let $\O_\q$ be the completion of $\O_K$ at $q$. We have
$$\lim_{n\ra\infty} \frac{|\left(\prod_{\q\mid \p} \O_\q^\ti\right)/\p^{n}|}{|\p^{nr}|} = \prod_{\q\mid\p} (1-\frac{1}{N_K(\q)}).$$
\end{lemma}

\begin{proof} 
Since any $\p$ is regular in $\O_K$, we have $\p = \prod_{\q\mid \p} \q^{e(\q)}$, where $e(\q)$ is the ramification index of $\q$. This gives
\begin{align*}
&\lim_{n\ra\infty} \frac{|\left(\left(\prod_{\q\mid\p} \O_\q \right)/\p^n\right)^\ti |}{|\p|^{nr}} \\
&= \lim_{n\ra\infty} \frac{|\left(\left(\prod_{\q\mid\p} \O_\q\right)/\prod_{\q\mid\p} \q^{e(\q)n} \right)^\ti|}{|\p|^{nr}} \\
&\stackrel{(1)}{=} \lim_{n\ra\infty} \frac{|\prod_{q\mid\p} (\O_\q/\q^{e(\q)n})^\ti |}{(N_F(\p))^{nr}} \\
&\stackrel{(2)}{=} \lim_{n\ra\infty} \frac{\prod_{\q\mid\p} ((N_K(\q))^{e(\q)n}-(N_K(\q))^{e(\q)n-1}) }{(N_F(\p))^{nr}} \\
&= \lim_{n\ra\infty} \frac{\prod_{\q\mid\p}\frac{N_K(\q)-1}{N_K(\q)} (N_K(\q))^{e(\q)n}}{(N_F(\p))^{nr}} \\
&= \lim_{n\ra\infty} \left(\prod_{\q\mid \p} 1-\frac{1}{N_K(\q)}\right) \frac{(N_K(\p))^n}{(N_F(\p))^{nr}} \\
&= \prod_{\q\mid\p} (1- \frac{1}{N_K(\q)}) \lim_{n\ra\infty} \frac{(N_F(\p))^{[K:F]n}}{(N_F(\p))^{rn}} \\
&= \prod_{\q\mid\p}(1-\frac{1}{N_K(\q)}). 
\end{align*}
Here $(1)$ holds, since $(a_\q)_{\q\mid\p}$ is a unit in $\prod_{\q\mid\p} \O_\q$ if and only if all $a_\q$ are units in $\O_\q$. Furthermore $(2)$ holds since $a \in \O_\q/\q^{e(\q)n}$ is not a unit if $\q \mid a$, which happens $N_K(\q)^{e(\q)n-1}$ times.  
\end{proof} 

\begin{prop} \th\label{LFcompfact} 
Let $f$ be irreducible. We have 
$$v_\p(f) = m_\p \cdot \frac{1-\frac{1}{|\p|}}{\prod_{\q\mid \p}(1-\frac{1}{N(\q)})},$$
where $m_\p$ is the number of $\GL_r(A_\p)$-orbits in the set $\{M \in \Mat_r(A_\p)\mid \charpoly(M)=f\}$.  
\end{prop}

\begin{proof}
Let $[M_0],\ldots,[M_{m_\p}]$ denote the $\GL_r(A_\p)$-orbits in the set $\{M \in \Mat_r(A_\p)\mid \charpoly(M)=f\}$. We denote by $[M_i](A/\p^n)$ the projection of $M_i$ onto $A/\p^n$. Note that this is an orbit under the action of $\GL_r(A/\p^n)$. Since orbits are disjoint, we get for all $n\geq 1$ that $|\{M\in \Mat_r(A/\p^n)\mid \charpoly(M)=f\}| = \sum_{i=0}^{m_\p} |[M_i](A/\p^n)|$. For each $i$ we have by the orbit-stabilizer theorem the exact sequence
$$1 \ra \Stab_{\GL_r(A/\p^n)}(M_i) \ra \GL_r(A/\p^n) \ra [M_i](A/\p^n) \ra 1.$$
This gives us that 
$$|\Stab_{\GL_r(A/\p^n)}(M_i)|\cdot |[M_i](A/\p^n)| = |\GL_r(A/\p^n)|.$$ 
This implies 
$$\sum_{i=0}^{m_\p} |[M_i(A/\p^n)])| = \sum_{i=0}^{m_\p} \frac{ |\GL_r(A/\p^n)|}{|\Stab_{\GL_r(A/\p^n)}(M_i)|}.$$
With \th\ref{CSstabM1} there exists some $N_i \in \GL_r(F_\p)$ such that $\Stab_{\GL_r(A_\p)}(M_i) = N_i\Stab_{\GL_r(A_\p)}(M_0)N_i^{-1}$. Since conjugation does not change the number of points, we get $|\Stab_{\GL_r(A/\p^n)}(M_i)| = |\Stab_{\GL_r(A/\p^n)}(M_0)|$ for all $i$. \\
\th\ref{CSstabis} gives that $\Stab_{\GL_r(A_\p)}(M_0) = \left(\prod_{\q\mid\p}\O_\q\right)^\ti$. With \th\ref{MprodOq} we get $$\lim_{n\ra\infty}\frac{\Big\vert \left(\left(\prod_{\q\mid\p}\O_\q \right)/\p^n\right)^\ti\Big\vert}{|\p|^{nr}} =  \prod_{\q\mid\p} \left(1-\frac{1}{N(\q)}\right).$$
Rewriting the result of \th\ref{MGLandSL} gives 
$$\lim_{n\ra\infty} \frac{|\GL_r(A/\p^n)|}{|SL_r(A/\p^n)|\cdot|\p|^{n}} = 1-\frac{1}{|\p|}.$$ 
Combining everything gives 
\begin{align*}
	v_\p(f) &= \lim_{n\ra\infty} \frac{|\{M \in \Mat_r(A/\p^n)\mid \charpoly(M)=f\}|\cdot |\p|^{(r-1)n}}{|\SL_r(A/\p^n)|} \\
	&= \lim_{n\ra\infty} \frac{ m_\p\cdot |\GL_r(A/\p^n)|\cdot |\p|^{rn}}{|\Stab_{\GL_r(A/\p^n)}|\cdot |\SL_r(A/\p^n)|\cdot |\p|^{n}} \\
	&= m_\p \frac{ 1-\frac{1}{|\p|}}{\prod_{\q\mid\p}(1-\frac{1}{N(\q)})}.
\end{align*}
\end{proof}

\begin{thm} \th\label{LFcompprod} 
Let $f$ be irreducible. We have an algorithm to compute $\prod_{\p} v_\p(f)$.
\end{thm}

\begin{proof}
With \th\ref{LFcompfact} we get that 
$$\prod_{\p} v_\p(f) = \prod_{\p} m_\p \frac{1-\frac{1}{|\p|}}{\prod_{\q\mid\p} 1-\frac{1}{N_\q}}.$$ 
With \th\ref{LFcomporbits} we can compute the number of orbits $m_\p$ for all primes $\p$ of $A$. Since for all but finitely many $\p$ we have $m_\p = 1$, we can compute their product $\prod_{\p}m_\p$. \\
Let $\F_{q^m}$ be the full constant field of $K$, so $\F_{q^m} = K \cap \ol{\F_q}$. Then, since every prime $\q$ of $K$ lies above one prime $\p$ of $A$, we get
$$\prod_{\p} \frac{1-\frac{1}{|\p|}}{\prod_{\q\mid\p} 1-\frac{1}{N(\q)}} = \frac{\prod_{\p} 1-\frac{1}{|\p|}}{\prod_{\q} 1-\frac{1}{N_K(\q)}} = \frac{\zeta_K(s)}{\zeta_F(s)}\Big\vert_{1}.$$ 
From \cite[Theorem 5.9]{Rosen} we know that there exist, explicitly computable polynomials $L_K$ and $L_F$ such that $\zeta_K(s) = \frac{L_K(q^{-ms})}{(1-q^{-ms})(1-q^{1-ms})}$ and $\zeta_F(s) = \frac{L_F(q^{-s})}{(1-q^{-s})(1-q^{1-s})}$. 
This gives us
$$\frac{\zeta_K(q^{-ms})}{\zeta_F(q^{-s})}\Big\vert_1 = \frac{L_K(q^{-ms})(1-q^{-s})(1-q^{1-s}) }{L_F(q^{-s})(1-q^{-ms})(1-q^{m(1-s)})}\Big\vert_1 = \frac{L_K(q^{-m})}{L_F(q^{-1})} \cdot \frac{1-q^{-1}}{1-q^{-m}} \cdot \frac{1}{m}.$$
\end{proof}

\section{The Characteristic 0 case} \label{secChar0} 
In this section we describe what happens if we change to the characteristic 0 setting. For the part about the Gekeler ratios, this means the objects we count are elliptic curves or abelian varieties of higher dimensions. For elliptic curves Gekeler already studied the isogeny classes and related those to Gekeler ratios, hence the name. 

\subsection{The Ideal Class Monoid over Number Fields}
In this setting $f(x)$ is an irreducible element of $\Z[x]$, $p \in \Z$ is a prime number, $K$ is the fraction field of $\Z[x]/f(x)$ and $R$ is a $\Z$-order in $K$. All results in Section \ref{secICM} carry over. We will state the most important ones here.

\begin{prop}[Analogue of \th\ref{spPicTrivRp}] The Picard group of $R_{\p}$ is trivial.
\end{prop}

\begin{prop}[Analogue of \th\ref{OObijApandAorders}]\th\label{Char0overorders} Let $f \in \Z[x]$ be irreducible and let $R = A[x]/f(x)$. Then the map
	$$\phi: \{p\text{-overorder of }R\} \ra \{\Z_{(p)}\text{-overorders of }R_{(p)}\},\ S \ra S_{(p)}$$
is a bijection.
\end{prop}

\begin{lemma}[Analogue of \th\ref{WEcalcLocWE}] \th\label{char0calcLocWE} 
	Let $I$ and $J$ be fractional $R$-ideals. Then $I_{(\p)} \sim J_{(\p)}$ if and only if 
	$$((I:J)\cdot (J:I)) \cap (R \bs [\p_1\cup \ldots \cup \p_k]) \not= \emptyset.$$
\end{lemma}

\begin{lemma}[Analogue of \th\ref{WEcompsurjonWES}] \th\label{char0compsurjonWES} 
	Suppose that $f$ is irreducible and that $R = \Z[x]/f(x)$. Let $S$ be the $p$-overorder of $R$ corresponding to $S_p$ under the map of \th\ref{Char0overorders}. Then the map $\theta: W_S(R) \ra W_{S_p}(R_p)$ defined by $[I] \ra [I_p]$ is well-defined and surjective.
\end{lemma}

In particular we still have: 
\begin{thm}[Analogue of \th\ref{ICMcomputable}] \th\label{char0ICMcomputable} 
	Let $f \in \Z[x]$ be irreducible and let $R = \Z[x]/f$. We can compute $ICM(R_\p)$. 
\end{thm} 

\begin{remark}
\th\ref{char0compsurjonWES,char0compsurjonWES} give an algorithm to compute $W_{S_p}(R_p)$. Together with \th\ref{Char0overorders} this means we can compute $W(R_p)$. In a different way, \cite[Lemmas 3.4 and 3.5 and Propositions 3.9 and 3.10]{AValgorithmsBKM}  also gives an algorithm to compute this. Let $\q_1,\ldots,\q_l$ be the prime ideals of $R$ that are both above $p$ and above the conductor $(R:\O_K)$ and let $n$ be the $p$-adic valuation of the index of $R$ in $\O_K$. They prove that there is a bijection between $W(R_p)$ and $\prod_{i=1}^l W(R+\q_i^{n}\O_K)$.\\
\th\ref{ICMisWE} tells us that $\ICM(R_p) = W(R_p)$, so this also provides an algorithm to compute the ideal class monoid without computing the local overorders. This means it also works if $R$ is not of the form $\Z[x]/f(x)$ with $f$ irreducible. \\
This algorithm requires to compute the weak equivalence classes of several orders, while in our algorithm we only need to compute (a subset of) the weak equivalence classes of $R$. This means we interchange computing overorders for computing more weak equivalence classes.     
\end{remark}

\subsection{Gekeler Ratios} 

Denote the weighted seize of an isogeny class of elliptic curves over $\F_p$ with trace of Frobenius $\tr(E) = a$ by 
$$H^*(a,p) = \sum_{E/\F_p,\ \tr(E) = a} \frac{1}{|\Aut_{\F_p}(E)|}.$$ In \cite[Theorem 5.5]{GekelerEllipticCurves} Gekeler proves that for an ordinary elliptic curve we have\\
$$H^*(a,p) = \frac{\sqrt{p}}{2} v_{\infty}(a,p) \prod_{l \text{ prime}} v_l(a,p).$$
Here $v_\infty(a,p) = \frac{2}{\pi} \sqrt{1-\frac{a^2}{4p}}$ is the global term and $$v_l(a,p) = \lim_{n \ra \infty} \frac{|\{M \in \Mat_2(\Z/l^k\Z) \mid \charpoly(M) = x^2-ax+p\}|}{l^{2k-2}(l^2-1)}$$ is the local Gekeler ratio at the prime $l$. \\
Write $D = a^2-4p$ and write $\d_2 = \d_2(a,p) = \max_{i \in \Z}(2^{2i}\mid D\text{ and }D\equiv 0,1 \mod 4)$ and $\d_l = \d_l(a,p) = \max_{i \in \Z}(l^{2i}\mid D)$. Then \cite[Corrolary 4.6]{GekelerEllipticCurves} gives an explicit formula for the local factors:
$$v_l(a,p) = (1-l^2)^{-1} \cdot \left(1 + l^{-1} + \begin{cases} 0 &\text{if } \qa{D/l^{2\d_l}}{l} = 1 \\ -(l+1)l^{-\d_l-2} &\text{if } \qa{D/l^{2\d_l}}{l} = 0 \\ -2l^{-\d_l-1} &\text{if } \qa{D/l^{2\d_l}}{l} = -1 \\ \end{cases} \right) .$$ 
The fundamental discriminant $D_0$ is defined by $D = c^2D_0$, where $c$ is maximal such that $c^2|D$ and $\frac{D}{c^2} \equiv 0,1 \mod 4$. In the proof of \cite[Corrolary 4.8]{GekelerEllipticCurves} Gekeler proves that if $l\nmid D$, then $v_l(a,p) = (1-\X(l)l^{-1})^{-1}$, where $\X$ is the character from the Kronecker symbol associated to $D_0$. This is a Dirichlet character by \cite[Theorem 9.13]{MultNTMontgomeryVaughan}. Furthermore we get by \cite[Section 9.5]{MultNTMontgomeryVaughan} that $L(s,\X_{D_0}) = \frac{\zeta_K(s)}{\zeta(s)}$, where $K$ is the quadratic field associated to $D_0$, i.e. $K = \Q(\sqrt{D_0})$. The residue of $\zeta(s)$ at $1$ is $1$ and the residue of $\zeta_K(s)$ can, assuming the generalized Riemann hypothesis, be computed up to arbitrary precision using the algorithms of \cite{ZetaFunction}. Since the size of the automorphism group of an elliptic curve is bounded, this is enough to get a precise answer. \\
This means that although our proof carries over to the new setting, we do not get anything new for elliptic curves. For completeness we do explain the differences. \\
Since in characteristic $0$ the irreducible $f$ is separable, we get that $f \not= g^2h$ for $g,h \in \Q_p[x]$. This gives the characteristic 0 version of \th\ref{CSfnotramified}. Write $f = f_1\cdots f_k$ with $f_i \in \Z_p[x]$ irreducible. Then the $f_i$ are distinct, since $f$ is separable. This means that $K_p = \prod_{\q | p} K_\q$, where $\q$ is a prime above $p$ as before, so \th\ref{CSsetisintcl} also holds in the characteristic 0 case. Together this gives:
\begin{thm}
Let $f$ be irreducible. We have 
$$v_p(f) = m_p \cdot \frac{1-\frac{1}{|p|}}{\prod_{\q\mid p}(1-\frac{1}{N(\q)})},$$
where $m_p$ is the number of $\GL_r(\Z_p)$-orbits in $\{M \in \Mat_r(\Z_p)\mid \charpoly(M)=f\}$ and $\q$ runs over the primes above $p$ in $K = \Frac(\Z[x]/f(x))$.\\
Furthermore
$$\prod_{p} v_p(f) = \frac{\zeta_K(s)}{\zeta(s)}\Big\vert_1 \prod_{p\text{ sing}} m_p,$$
where the second product is over all primes that are singular in $\Z[x]/f(x)$. 
\end{thm}
As explained above we can compute the residues up to finite precision. Since the global term is also explicit, this is enough to compute the number of isomorphism classes of elliptic curves inside an isogeny class. \\
For abelian varieties of higher dimension, Achter, Altuğ, Garcia and Gordon have in \cite{GekelerRatiosAV} a similar result. Let $X$ be a principally polarized abelian variety with commutative endomorphism ring. If $X$ is ordinary or $X$ is defined over $\F_p$, then 
$$H^*(X,\F_q) = q^{\frac{1}{2}g(g+1)}\tau_Tv_\infty(X,\F_q) \prod_{l\text{ prime}} v_l(X,\F_q),$$
where $g = \dim(X)$, $\t_T$ is a Tamagawa number of a torus $T$ associated to $X$, $v_\infty(X,\F_q)$ is a global term and $v_l(X,\F_q)$ are the local factors. Let $f$ be the characteristic polynomial of the Frobenius of $X$. If $l\nmid p\cdot \text{disc}(f)$, then
$$v_l(X,\F_q) = \lim_{n\ra \infty} \frac{ |\{M \in \GSp_{2g}(\Z_l/l^n) \mid \charpoly(M) = f\}|}{|\GSp_{2g}(\Z_l/l^{n})|\ /\ |\mathbb{A}_{\GSp_{2g}}(\Z/l^n)|}.$$
Here $\GSp_{2g}(S) = \{M \in \Mat_{2g}(S) \mid M^t J M = J\}$ with $J = \mat{0& -I_g\\I_g & 0}$ and $\mathbb{A}_{\GSp_{2g}}$ is the space of all possible characteristic polynomials of matrices in $\GSp_{2g}$. \\ 
Since we need to count over $\GSp_{2g}$ instead of over all matrices we cannot use the results from the function field setting. The first obstruction to get the method to work is that there is no description of the stabilizers of $\GSp_{2g}$ acting on itself that allows for the point counting. If one has such a description, this can maybe be used to formulate and prove an analogue of \th\ref{LFWeiXu} for $\GSp_{2g}$. 

\section{Implementation} \label{secImplementation} 

The code can be find in \cite{codeLocalGekelerRatios} and runs in Magma, \cite{Magma}. The code implements the algorithms described in this paper in the function field setting. Since neither Magma nor SageMath can handle inseparable function field extensions properly, this case is not yet implemented. The implementation builds on code of Katen, see \cite{codeJeffKaten}.\\
The implementation we use to compute all $\p$-overorders generates all submodules of $\O/R$ where $\O$ is the maximal $\p$-overorder and checks if the pullback of those is an order. In \cite{OverordersHS} a faster algorithm is described which generates the overorders recursively. This can be used to get a faster implementation, but is not yet implemented in the function field setting. \\
In the separable case there is also a new algorithm, see \cite[Algorithm ComputeWbar]{StefanoWeakEquiv} to compute $W_S(R)$. In the characteristic $0$ setting this algorithm performed better in practice then the previous one. Since also this algorithm is not implemented in the function field setting, we use the algorithm described in \cite{Stefano}, which was already implemented by Katen.

	\bibliographystyle{amsplain} 
	\bibliography{biblio}
	
\end{document}